\documentclass[final, 3p, times]{elsarticle}%
\usepackage{graphicx}
\usepackage{xcolor}
\usepackage{colortbl}
\usepackage{epstopdf}
\usepackage{amssymb}
\usepackage{amsthm}
\usepackage{amsmath}
\usepackage{color}
\usepackage{dsfont}
\usepackage{mathtools}
\usepackage{extarrows}
\usepackage{hyperref}
\usepackage{bm}
\usepackage{caption}
\usepackage{float}
\usepackage{verbatim}
\usepackage{epsfig}
\usepackage{multirow}
\usepackage{subfig}
\usepackage[section]{placeins}
\usepackage{float}
\usepackage{cleveref}
\usepackage{xr}

\usepackage{algorithm}       % 定义 algorithm 浮动体环境
\usepackage{algpseudocode}   % 提供 algorithmc 代码排版（或 algorithmicx）
\newtheorem{exm}{Example}[section]

\biboptions{compress}

\journal{Elsevier}

\begin{document}
\begin{frontmatter}
\title{ A note on spectral Monte-Carlo method for fractional Poisson equation on high-dimensional ball}

\tnotetext[label1]{The research of Dongling Wang is supported in part by the National Natural Science Foundation of China under grants 12271463. The work of Lisen Ding is supported by Postgraduate Scientific Research Innovation Project of Xiangtan University, China (No. XDCX2025Y213). The work of Mingyi Wang is supported by Postgraduate Scientific Research Innovation Project of Hunan Province, China (No. CX20240602).
 \\ Declarations of interest: none.
}

\author[XTU]{Lisen Ding}
\ead{dingmath15@smail.xtu.edu.cn}
\author[XTU]{Mingyi Wang}
\ead{202331510144@smail.xtu.edu.cn}
\author[XTU]{Dongling Wang\corref{mycorrespondingauthor}}
\ead{wdymath@xtu.edu.cn; ORCID: 0000-0001-8509-2837}
\cortext[mycorrespondingauthor]{Corresponding author.}

\address[XTU]{Hunan Key Laboratory for Computation and Simulation in Science and Engineering, School of Mathematics and Computational Science, Xiangtan University, Xiangtan, Hunan 411105, China}

\begin{abstract}

Recently, a class of efficient spectral Monte Carlo methods was developed in \cite{Feng2025ExponentiallyAS} for solving fractional Poisson equations. These methods fully consider the low regularity of the solution near boundaries and leverage the efficiency of walk-on-spheres algorithms, achieving spectral accuracy. However, the underlying formulation is essentially one-dimensional.
In this work, we extend this approach to radial solutions in general high-dimensional balls. This is accomplished by employing a different set of eigenfunctions for the fractional Laplacian and deriving new interpolation formulas. We provide a comprehensive description of our methodology and a detailed comparison with existing techniques. Numerical experiments confirm the efficacy of the proposed extension.

\end{abstract}


\end{frontmatter}

\section{Introduction}

We consider the fractional Poisson equation with exterior Dirichlet boundary conditions:
\begin{equation}\label{eq:fracPoiss}
\begin{cases}
(-\Delta)^{\frac{s}{2}} u(x) = f(x), & x \in \Omega, \\
u(x) = g(x), & x \in \Omega^c,
\end{cases}
\end{equation}
where $\Omega \subset \mathbb{R}^n$ is a bounded domain and $f: \Omega \to \mathbb{R}$ is a given function. The integral fractional Laplacian for $s \in (0,2)$ is defined as:
\begin{equation}\label{eq:fracLap}
(-\Delta)^{\frac{s}{2}} u(x) = C_{n,s} , \mathrm{P.V.} \int_{\mathbb{R}^n} \frac{u(x) - u(y)}{|x - y|^{n + s}} , \mathrm{d}y, \quad C_{n, s} = \frac{2^s \Gamma\left(\frac{n + s}{2}\right)}{\pi^{\frac{n}{2}} \Gamma\left(1 - \frac{s}{2}\right)},
\end{equation}
where P.V. denotes the Cauchy principal value. This model has garnered significant attention for its theoretical and practical relevance \cite{fernandez2024integro, d2020numerical, bonito2018numerical}.
A key feature of this operator is the singularity of its kernel, which induces low solution regularity near the boundary regardless of the smoothness of $f$ or $\partial \Omega$. This inherent singularity, combined with the nonlocal nature of the operator, poses fundamental challenges for numerical computation, particularly in high dimensions.

Conventional numerical methods, such as finite difference \cite{huang2014numerical, hao2021fractional, minden2020simple, DUO2019639, huang2024grid, han2022monotone} and finite element \cite{Acosta2018Regular, bonito2018numerical, bahr2024implementation} methods, often lead to complex discretization and dense matrices. The standard spectral method \cite{hao2021sharp, xu2018spectral} also faces significant challenges in high-dimensional settings.

To address these issues, a walk-on-spheres (WOS) algorithm was introduced in \cite{Kyprianou2017Unbiased} for two-dimensional problems. This method leverages a probabilistic representation of the solution, which is then computed via Monte-Carlo simulation. Its key advantage lies in handling complex geometries effectively. This approach has since been extended to higher dimensions ($n \geq 3$) and more complex semilinear models \cite{Sheng2023Efficient, sheng2024implicit, Jiao2023Modified, yuan2025class, nie2025walk}. However, the Monte Carlo nature of WOS typically limits its convergence rate to half-order accuracy.

A recent breakthrough in \cite{Feng2025ExponentiallyAS} achieved spectral accuracy by constructing sophisticated interpolation operators that exploit the precise Hölder regularity of the solution near the boundary. Nevertheless, this method is fundamentally one-dimensional in its construction (see Section \ref{sec:1eig}). In this article, we generalize this high-accuracy approach to compute radial solutions in high-dimensional spherical domains.

The paper is organized as follows. Section \ref{sec:pre} discusses radial solutions for problem \eqref{eq:fracPoiss} and outlines the relevant stochastic algorithms. Section \ref{sec:jacobi} provides necessary background on Jacobi polynomials. The one-dimensional interpolation formula from \cite{Feng2025ExponentiallyAS} and our new interpolation formula for radial solutions on the unit ball are presented in Sections \ref{sec:1eig} and \ref{sec:neig}, respectively. A comparative analysis and discussion of these two approaches are given in Section \ref{sec:Compa}. Numerical experiments in Section \ref{sec:numexp} demonstrate the effectiveness of our method. Finally, concluding remarks are offered in Section \ref{sec:conclude}.

\section{Preliminaries}\label{sec:pre}
In what follows, we define $B_r^{x_0}:= \{x \in \mathbb{R}^n : |x - x_0| \leq r\}$ for any $x_0 \in \mathbb{R}^n$, where $|\cdot|$ denotes the standard Euclidean norm in $\mathbb{R}^n$.

\subsection{Radial solutions}

Radial solutions serve as an important tool in the study of partial differential equations (PDEs), often enabling explicit expressions of solutions. For equations with complex structure, radial solutions offer a simplified approach, and their analysis can lead to deeper insights into the behavior of such equations. The fractional Laplacian \eqref{eq:fracLap}, defined as a hypersingular integral operator, has been extensively studied in the context of radial solutions; see, for example, \cite{felmer2014radial, Chen2021move, Frank2016Uniqueness, Dyda2012power}. In what follows, we discuss some key aspects of radial solutions related to the integral fractional Laplacian.

It is well-known that the fractional Laplacian \eqref{eq:fracLap} reduces to the classical negative Laplacian $-\Delta$ when $s=2$. As noted in \cite{evans2022partial}, the rotational invariance of $-\Delta$ implies that the equation $-\Delta u = 0$ in $\mathbb{R}^{n}$ admits radial solutions. Specifically, if $R \in \mathbb{R}^{n\times n}$ is an orthogonal matrix and $v(x) = u(Rx)$, then $-\Delta v = 0$.
This property extends to the fractional operator $(-\Delta)^{\frac{s}{2}}$, which is also rotationally invariant due to the dependence of its kernel $|x-y|^{-(n+s)}$ solely on the distance between points. To verify this, suppose $(-\Delta)^{\frac{s}{2}}u = 0$ and define $v(x) = u(Rx)$. Then,
\begin{align*}
(-\Delta)^{\frac{s}{2}}v(x) &= C_{n,s} \lim_{\varepsilon \to 0} \int_{\mathbb{R}^n \setminus B_\varepsilon(x)} \frac{v(x) - v(y)}{|x - y|^{n+s}}\mathrm{d}y= C_{n,s} \lim_{\varepsilon \to 0} \int_{\mathbb{R}^n \setminus B_\varepsilon(x)} \frac{u(Rx) - u(Ry)}{|x - y|^{n+s}}\mathrm{d}y.
\end{align*}
Making the change of variables $z = Ry$ and using the orthogonality of $R$ (which gives $|x-y| = |Rx - z|$ and $\mathrm{d}y = \mathrm{d}z$), we obtain
\begin{align*}
(-\Delta)^{\frac{s}{2}}v(x) &= C_{n,s} \lim_{\varepsilon \to 0} \int_{\mathbb{R}^n \setminus B_\varepsilon(Rx)} \frac{u(Rx) - u(z)}{|Rx - z|^{n+s}} \mathrm{d}z= (-\Delta)^{\frac{s}{2}}u(Rx) = 0.
\end{align*}
Hence, $(-\Delta)^{\frac{s}{2}} v = 0$, confirming that the equation is rotationally invariant and therefore admits radial solutions. Explicit forms of such radial solutions can be found in \cite{chen2020fractional}.

However, for a bounded domain $\Omega$ with a non-trivial source term $f \neq 0$, the existence of a radial solution to problem \eqref{eq:fracPoiss} requires both $\Omega$ and $f$ to satisfy specific symmetry conditions. A common framework, widely adopted in the literature and in this work, is to impose the homogeneous Dirichlet condition $g=0$.
Leveraging the rotational invariance of $(-\Delta)^{\frac{s}{2}}$, if $\Omega$ is a ball (e.g., the unit ball $B_1^0$) or other spherical domain and $f$ is a radial function, then the solution to \eqref{eq:fracPoiss} is also radial \cite{felmer2014radial, Dyda2012power}. Since this article is restricted to a ball, it suffices to assume that $f$ is radial to guarantee a radial solution.

\subsection{Stochastic algorithm for fractional Possion problem}

For spherical domains $\Omega = B_r^0$, the solution to problem \eqref{eq:fracPoiss} admits the representation
\begin{equation}\label{eq:ballsolution}
u(x) = \int_{B_r^{0}} f(y) Q_r(x,y) \mathrm{d}y + \int_{\mathbb{R}^n \setminus B_r^{0}} g(z) P_r(x,z) \mathrm{d}z, \quad x \in B_r^{0},
\end{equation}
where $Q_r$ is the fractional Green's function and $P_r$ is the fractional Poisson kernel. As shown in \cite{Sheng2023Efficient}, this representation can be reformulated in terms of expectations:
\begin{equation}
u(x) = \zeta(x) \mathbb{E}_{\widetilde{Q}r}[f(Y)] + \mathbb{E}_{P_r}[g(Z)], \quad Y \in B_r^{0},\ Z \in \mathbb{R}^n \setminus B_r^{0},
\end{equation}
where the weight function $\zeta(x) = \int{B_r^{0}} Q_r(x,y) \mathrm{d}y$ for $x \in B_r^{0}$, and $ \widetilde{Q}r(x,y) = Q_r(x,y)/\zeta(x)$. Combining this expectation formula with the walk-on-spheres algorithm, the solution on a general domain can be expressed as a conditional expectation. The corresponding Monte-Carlo discretization yields the numerical scheme:
\begin{equation}\label{eq:alma}
u{M}(x) = \frac{1}{M}\sum_{i=1}^{M}S^{i} = \frac{1}{M}\sum_{i=1}^{M}\left[\sum_{k=0}^{L-1} \zeta(X_{k}^{i})f(Y_{k+1}^{i}) +g(X_{L}^{i})\right],
\end{equation}
where $S^i$ denotes the $i$-th sample path in a total of $M$ Monte Carlo trials. For each sample, the walk-on-spheres algorithm generates a sequence of balls $\big\{B^{X_k}_{r_{k+1}}\big\}_{k=0}^{L-1}$, with radii $r_{k+1} = \mathrm{dist}(X_k,\partial\Omega)$. Here, $X_k$ represents the position of the Lévy process $X_{\tau}^s$ at time $\tau_k$, and $Y_{k+1}$ is a point sampled inside the ball $B^{X_k}_{r{k+1}}$. The algorithm terminates after a finite number of steps $L$, ensuring its practical feasibility.

A variant of this scheme is presented in \cite{Jiao2023Modified}, which differs from \eqref{eq:alma} in the choice of the weight function and the treatment of the source term. Defining
\begin{equation*}
\tilde{f}(Y_{k+1}) = \left[1 - B\left(\frac{\gamma^2}{r_k^2}; \frac{n-s}{2},1-\frac{s}{2}\right)\right]f(Y_{k+1}), \quad
\eta(X_k) = \frac{B\left(\frac{n-s}{2},\frac{s}{2}\right)}{2^{s-1}\Gamma^2(s/2)}r_{k+1}^s,
\end{equation*}
where $B(a,b)$ is the beta function and $B(\cdot; a,b)$ is the incomplete beta function with $B(1; a,b)=1$, the numerical formula in \cite{Jiao2023Modified} is obtained by replacing $\zeta(X_{k})$ and $f(Y_{k+1})$ in \eqref{eq:alma} with $\eta(X_k)$ and $\tilde{f}(Y_{k+1})$, respectively. A key advantage of this modified scheme is that the weight $\eta(X_k)$ is given in closed form, eliminating the need for numerical integration. For this reason, we adopt the formulation from \cite{Jiao2023Modified} in the present work. Further details on the stochastic algorithms for \eqref{eq:fracPoiss} can be found in \cite{Sheng2023Efficient, Jiao2023Modified}.

\section{Spectral Monte Carlo Method}\label{sec:SMC}

In this section, we always employ $t\in\mathbb{R}$ as one-dimensional variables, and $x\in\mathbb{R}^n$ as high-dimensional variables.
We recall some basic notation and properties of Jacobi polynomials \cite{shen2011spectral}. 

\subsection{Jacobi Polynomials}\label{sec:jacobi}

The Jacobi polynomials $\left\{ P_m^{\alpha,\, \beta}(t) \right\}$ with $\alpha,\beta > -1$ are mutually orthogonal with respect to the weight function $\omega^{\alpha,\,\beta}(t) = (1-t)^\alpha(1+t)^\beta$:
\begin{equation}\label{eq:Jac1}
  \int_{-1}^1 \omega^{\alpha,\,\beta}(t) P_m^{\alpha,\,\beta}(t) P_q^{\alpha,\,\beta}(t) \ \mathrm{d}t = \gamma_m^{\alpha, \,\beta}\delta_{mq},
\end{equation}
where $\delta_{mq}$ denotes the Kronecker function, and
$
\gamma_m^{\alpha, \, \beta} = \frac{2^{\alpha+\beta+1}\Gamma(m+\alpha+1)\Gamma(m+\beta+1)}{(2m+\alpha+\beta+1)\Gamma(m+\alpha+\beta+1)m!}.\
$
 These polynomials can be recursively defined by:
\begin{equation*}
  \begin{cases}
    P_0^{\alpha,\,\beta}(t) = 1, \\
    P_1^{\alpha,\,\beta}(t) = \frac{\alpha+\beta+2}{2}t + \frac{\alpha-\beta}{2}, \\
    P_m^{\alpha,\,\beta}(t) = (A_m t + B_m)P_{m-1}^{\alpha,\,\beta}(t) - C_m P_{m-2}^{\alpha,\,\beta}(t) \text{ for } m \geq 2,
  \end{cases}
\end{equation*}
where the coefficients 
$
    A_m  = \frac{(2m+\alpha+\beta-1)(2m+\alpha+\beta)}{2m(m+\alpha+\beta)}, 
    B_m  = \frac{(\alpha^2-\beta^2)(2m+\alpha+\beta-1)}{2m(m+\alpha+\beta)(2m+\alpha+\beta-2)}, 
    C_m  = \frac{(m+\alpha-1)(m+\beta-1)(2m+\alpha+\beta)}{m(m+\alpha+\beta)(2m+\alpha+\beta-2)}.
$
Leveraging the weighted orthogonality, we employ the Jacobi-Gauss quadrature. Let 
$\left\{t_k^{\alpha,\,\beta}, \omega_k^{\alpha,\,\beta} \right\}_{k=0}^{N_t}$ 
denote the Jacobi-Gauss quadrature nodes and weights,  which satisfy:
\begin{equation}\label{eq:20}
  \int_{-1}^1 \psi(t)\omega^{\alpha,\,\beta}(t) \mathrm{d}t = \sum_{k=0}^{N_t} \psi \left( t_k^{\alpha,\,\beta} \right)\omega_k^{\alpha,\,\beta}, \quad \forall \psi \in \mathcal{P}_{2N_t+1},
\end{equation}
where $\mathcal{P}_N$ represents the space of polynomials of degree up to $N$.

\subsection{Eigenfunctions and interpolation formula used in \cite{Feng2025ExponentiallyAS}}\label{sec:1eig}
When $n=1$ and $\Omega:=\Lambda=(-1,1)$, the eigenfunctions of the fractional Laplacian with respect to the weight function $\rho(t):=(1-t^2)^{\frac{s}{2}}$ are Jacobi polynomials $P_m^{\frac{s}{2},\,\frac{s}{2}}(t)$ according to \cite{MAO2016Efficient}, i.e.,
%\mg{
\begin{equation}\label{eq:eigfun1}
(-\Delta)^{\frac{s}{2}} \left[ \rho(t)P_m^{\frac{s}{2},\frac{s}{2}}(t) \right]  =\lambda_m P_m^{\frac{s}{2},\frac{s}{2}}(t), 
\quad \lambda_m = \frac{\Gamma(m+\alpha+1)}{m!},\quad  t\in \Lambda.
\end{equation}
This is the crucial starting point of spectral Monte-Carlo method stated in \cite{Feng2025ExponentiallyAS}. 
%}

Due to the low regularity of solutions to fractional Poisson equations near the boundary of the region, conventional interpolation operators are unable to achieve high-precision approximation. To overcome this difficulty, a class of generalized interpolation operators with fractional order factors are introduced in \cite{Feng2025ExponentiallyAS}. Using Jacobi-Gauss quadrature nodes $\left\{ t_k^{\alpha,\,\beta}, \omega_k^{\alpha,\,\beta} \right\}_{k=0}^{N_t}$ defined in \eqref{eq:20} with $\alpha=\beta=\frac{s}{2}$, one can define the generalized Lagrange interpolation operator:
\begin{equation}\label{eq:interp}
  \mathcal{I}_{N_t}^{\frac{s}{2}}v(t) = \sum_{k=0}^{N_t}v(t_k)l_k^{\frac{s}{2}}(t),
\end{equation}
with interpolating basis functions
\begin{equation}\label{eq:fract1}
  l_k^{\frac{s}{2}}(t) = \psi(t) \cdot h_k(t),\quad \psi(t) = \left(\frac{1-t^2}{1-t_k^2}\right)^{\frac{s}{2}},\
  h_k(t) = \prod_{\substack{0\leq i \leq N_t \\ i\neq k}} \frac{t-t_i}{t_k - t_i},
\end{equation}
where $h_k(t)$ are standard Lagrange basis polynomials. To simplify the notation, we have denoted $t_k := t_k^{\frac{s}{2},\,\frac{s}{2}}$. Next, we choose generalized Jacobi functions $J^{\frac{s}{2}}_m(t)=(1-t^2)^{\frac{s}{2}}P_m^{\frac{s}{2},\frac{s}{2}}(t)$ as the basis function, then $l_k^{\frac{s}{2}}(t)$ can be expressed as the following expansion by the Jacobi-Gauss quadrature formula \eqref{eq:20}, that is,
\begin{equation}\label{eq:1LagJac}
  l_k^{\frac{s}{2}}(t) = \sum_{m=0}^{N_t}c_{mk}J^{\frac{s}{2}}_m(t),\quad c_{mk} = \frac{1}{\gamma_m^{\frac{s}{2},\frac{s}{2}}}(1-t_k^2)^{-\frac{s}{2}}  P_m^{\frac{s}{2},\frac{s}{2}}(t_k) \omega_k^{\frac{s}{2},\,\frac{s}{2}}.
\end{equation}
Therefore, the interpolation function $\mathcal{I}_{N_t}^{\frac{s}{2}}v(t)$ can be rewritten as
\begin{equation}\label{eq:interp1}
  \mathcal{I}_{N_t}^{\frac{s}{2}}v(t) = \sum_{k=0}^{N_t}v(t_k)\left(\sum_{m=0}^{N_t}c_{mk}J^{\frac{s}{2}}_m(t)\right),
\end{equation}
and combine it with \eqref{eq:eigfun1}, we have
\begin{equation}\label{eq:interpfra}
  \begin{aligned}
    (-\Delta)^{\frac{s}{2}}\mathcal{I}_{N_t}^{\frac{s}{2}}v(t) &= \sum_{k=0}^{N_t}v(t_k)\left(\sum_{m=0}^{N_t}c_{mk}(-\Delta)^{\frac{s}{2}}J^{\frac{s}{2}}_m(t)\right) = \sum_{k=0}^{N_t}v(t_k)\left(\sum_{m=0}^{N_t}c_{mk} \lambda_m P^{\frac{s}{2},\frac{s}{2}}_m(t)\right) \\
    &=\sum_{m=0}^{N_t}\hat{v}_{m}P^{\frac{s}{2},\frac{s}{2}}_m(t), \quad \hat{v}_{m} = \lambda_m \sum_{k=0}^{N_t}c_{mk}  v(t_k).
  \end{aligned}
\end{equation}

Next, we will focus on the problem \eqref{eq:fracPoiss} with homogeneous boundary conditions $g(x)=0$. The spectral Monte-Carlo method derived from \eqref{eq:interp1} and \eqref{eq:interpfra} can be summarized as Algorithm \ref{alg:SMC}.

\begin{algorithm}
	\caption{Spectral Monte Carlo iteration algorithm}
	\label{alg:SMC}
	\begin{algorithmic}[1]
		\State \textbf{Input:} Jacobi-Gauss points $\{t_k\}_{k=0}^{N_t}$ and iteration number $K$.
    \For{$i=1$; $i\leq K$}
    \If{$i=1$}
    \State The initial solution $\left\{u^{(i)}_{\star}(t_k)\right\}_{k=0}^{N_t}$ can be calculated by stochastic formula \eqref{eq:alma}.
    \EndIf

    \State High accuracy approximation for $(-\Delta)^{\frac{s}{2}}u^{(i)}_{\star}(t)$ obtained by replacing $v(t_k)$ in \eqref{eq:interpfra} with $u ^{(i)} _ {\star} (t_k) $.

    \State \parbox[t]{\dimexpr\linewidth-\algorithmicindent}    {Set $\varepsilon=u-u ^ {(i)} _ {\star} $, then construct residual problem
    \begin{equation*}%\label{eq:29}
      \begin{cases}
        (-\Delta)^{\frac{s}{2}} \varepsilon = f^{(1)}, \quad &x \in \Omega, \\
        \varepsilon = 0, \quad &x \in \Omega^{c},
      \end{cases}
    \end{equation*}
    where $f ^ {(i)}=f - (- \Delta) ^ {\frac {s} {2}} u ^ {(i)} _ {\star} $. Again, $\left\{\varepsilon_{\star}^{(i)}(t_k)\right\} _ {k=0} ^ {N_t}$  can be evaluated by stochastic formula \eqref{eq:alma}.}

    \State Update the solution as $u ^ {(i+1)} _ {\star} (t)=u ^ {(i)} _ {\star} (t)+\varepsilon ^ {(i)} _ {\star} (t) $ by replacing $v(t_k)$ in \eqref{eq:interp1} with $u ^ {(i+1)} _ {\star} (t_k)$.
    
    \EndFor
		\State \textbf{Output:} $u^{(i+1)}_ {\star}(t)$.
	\end{algorithmic}
\end{algorithm}

\subsection{Eigenfunctions and interpolation formula of radial solutions}\label{sec:neig}

The algorithm presented in Algorithm \ref{alg:SMC} for solving problem \eqref{eq:fracPoiss} is inherently one-dimensional. This limitation arises from the definition of the integral fractional Laplacian. In full space $\mathbb{R}^n$, the operator \eqref{eq:fracLap} admits an equivalent definition via the Fourier transform:
\[
\mathcal{F}[(-\Delta)^{\frac{s}{2}}u](\xi) = |\xi|^s\mathcal{F}[u](\xi), \text{ where } \xi=(\xi_1,\xi_2,\cdots,\xi_n).
\]
When $s=2$, this reduces to the standard Fourier characterization of the classical Laplacian. While the identity $|\xi|^2 = \xi_1^2 + \xi_2^2 + \cdots + \xi_n^2$ allows the Laplacian to be naturally extended to $n$ dimensions, the fractional power $|\xi|^s$ does not permit a similar separable decomposition. Consequently, the spectral Monte-Carlo method based on the one-dimensional eigenfunction relation \eqref{eq:eigfun1} cannot be directly generalized to higher dimensions.

However, by treating $|\xi|$ as a single variable, we can instead consider the radial solutions of problem \eqref{eq:fracPoiss}. This perspective opens a pathway for extending the spectral Monte-Carlo method to higher-dimensional scenarios. To realize this, the one-dimensional eigenfunction relation \eqref{eq:eigfun1} must be replaced with a radial analogue. The derivation of such a formula is the focus of this section.

Let $\Omega := B_1^0 \subset \mathbb{R}^n$. In \cite{Dyda2017Eigenvalues}, the eigenfunctions of the fractional Laplacian with respect to the weight function $\rho(x):= (1-|x|^2)^{\frac{s}{2}}$ are Jacobi polynomials with indices $(\alpha, \beta)=(s/2, n/2-1)$, i.e.,
\begin{equation}\label{eq:eigfun}
(-\Delta)^{\frac{s}{2}} \left[ \rho(x)P_m^{\frac{s}{2},\frac{n}{2}-1}(2|x|^2-1) \right]  = \mu_m P_m^{\frac{s}{2},\frac{n}{2}-1}(2|x|^2-1),
\end{equation}
where $ \mu_m = \frac{2^s \Gamma\left(\frac{s}{2}+m+1\right)\Gamma\left(\frac{s+n}{2}+m\right)}{m!\Gamma\left(\frac{n}{2}+m\right)}.$
Using \eqref{eq:Jac1}, we can verify the orthogonality of Jacobi polynomials $P_m^{\frac{s}{2}, \frac{n}{2}-1}(2|x|^2-1)$ with weight $\rho(x)$:
  \begin{align}
    & \int_{B_1^0} \rho(x) P_m^{\frac{s}{2},\frac{n}{2}-1}(2|x|^2-1) P_q^{\frac{s}{2},\frac{n}{2}-1}(2|x|^2-1) \ \mathrm{d}x \notag\\
    &= \int_0^1 r^{n-1} (1-r^2)^{\frac{s}{2}} P_m^{\frac{s}{2},\frac{n}{2}-1}(2r^2-1)P_q^{\frac{s}{2},\frac{n}{2}-1}(2r^2-1) \ \mathrm{d}r \cdot \int_{\partial B_1^0} \mathrm{d}V \notag\\
    &=\frac{2\pi^{n/2}}{\Gamma\left(\frac{n}{2}\right)} \cdot 2^{-\frac{n+s}{2}-1} \int_{-1}^1 (1+t)^{n/2-1}(1-t)^{s/2}P_m^{\frac{s}{2},\frac{n}{2}-1}(t)P_q^{\frac{s}{2},\frac{n}{2}-1}(t) \ \mathrm{d}t \notag\\
    &= \frac{\pi^{n/2}}{2^{\frac{n+s}{2}}\Gamma\left(\frac{n}{2}\right)} \gamma_m^{\frac{s}{2},\frac{n}{2}-1}\delta_{mq}, \label{eq:orth}
  \end{align}
where $\gamma_m^{\frac{s}{2},\frac{n}{2}-1}$ is the Jacobi orthogonality constant from \eqref{eq:Jac1}.

Under the condition of problem \eqref{eq:fracPoiss} exists radial solutions, our goal is to make the interpolation operator \eqref{eq:interp} suitable for high-dimensional space with the help of variable relation $t=2|x|^2-1$ for $x\in B_1^0  \subset \mathbb{R}^n$. Based on this variable relation and using Jacobi-Gauss quadrature nodes $\left\{ t_k^{\alpha,\,\beta}, \omega_k^{\alpha,\,\beta} \right\}_{k=0}^{N_t}$ with index $\alpha=s/2$ and $\beta=n/2-1$, we construct the new factor term $\hat{\psi}(t) = \left(1-t\right)^{\frac{s}{2}}/\left(1-t_k\right)^{\frac{s}{2}}$ as displayed in \eqref{eq:fract1}, where $t_k:=t_k^{\frac{s}{2},\frac{n}{2}-1}$. This consideration is the difference between our approach and the one in \cite{Feng2025ExponentiallyAS}. For more detailed explanations and analyses, please refer to Section \ref{sec:Compa} below.

In order to use the orthogonality relation \eqref{eq:orth}, we firstly need to build the connection between Lagrange polynomials $h_k(t)$ and Jacobi polynomials $P_m^{\frac{s}{2},\frac{n}{2}-1}(t)$.  Expanding $h_k(t)$ as $h_k(t) = \sum_{m=0}^{N_t} \widehat{C}_{m,k} P_m^{\frac{s}{2},\frac{n}{2}-1}(t)$, 
the coefficients $\widehat{C}_{m,k}$ are determined via orthogonality:
\begin{equation}\label{eq:coe}
    \int_{-1}^1(1+t)^{\frac{n}{2}-1} (1-t)^{\frac{s}{2}} h_k(t) P_m^{\frac{s}{2},\frac{n}{2}-1}(t) \ \mathrm{d}t 
    = \widehat{C}_{m,k} \int_{-1}^1(1+t)^{\frac{n}{2}-1}(1-t)^{\frac{s}{2}}\left[P_m^{\frac{s}{2},\frac{n}{2}-1}(t)\right]^2 \mathrm{d}t 
    = \widehat{C}_{m,k} \gamma_m^{\frac{s}{2},\frac{n}{2}-1}.
\end{equation}
Applying Jacobi-Gauss quadrature \eqref{eq:20}, we have
\begin{equation*}
  \int_{-1}^1(1+t)^{\frac{n}{2}-1} (1-t)^{\frac{s}{2}} h_k(t) P_m^{\frac{s}{2},\frac{n}{2}-1}(t) \ \mathrm{d}t  = \sum_{i=0}^{N_t} \omega_i^{\frac{s}{2},\frac{n}{2}-1} h_k(t_i) P_m^{\frac{s}{2},\frac{n}{2}-1}(t_i).
\end{equation*}
Noting the nodal property $h_k(t_i) = \delta_{ik}$, we obtain that
$
\widehat{C}_{m,k} = \frac{1}{\gamma_m^{\frac{s}{2},\frac{n}{2}-1}} \omega_k^{\frac{s}{2},\frac{n}{2}-1} P_m^{\frac{s}{2},\frac{n}{2}-1}(t_k).
$

In addition, considering the low regularity of solutions to high-dimensional fractional Poisson equations near the boundaries, we need to factorize the factor term $\hat{\psi}(t)$ as 
\begin{equation}\label{eq:fract2}
\hat{\psi}(t) = \left(\frac{2}{1-t_k}\right)^{\frac{s}{2}}\left(\frac{1-t}{2}\right)^{\frac{s}{2}}.
\end{equation}
This factorization allows us to recover the weight  $\rho(x)=(1-|x|^2)^{\frac{s}{2}}$ by the variable relation $t=2|x|^2-1$ for $x\in B_1^0 \subset \mathbb{R}^n$.
It gives us the opportunity to include weight factors $\rho(x)$ when calculating the fractional Laplacian of the interpolation operator in \eqref{eq:interp} by the continuity equation \eqref{eq:eigfun}.
The above analysis inspired us to define a new interpolation operator
\begin{equation}\label{eq:ninterp}
\begin{aligned}
  \mathcal{I}_{N_t}^{\frac{s}{2}}v(t) &= \sum_{k=0}^{N_t}v(t_k)  \hat{\psi}(t)h_{k}(t) 
  =\sum_{k=0}^{N_t}v(t_k)  \left(\frac{2}{1-t_k}\right)^{\frac{s}{2}}\left(\frac{1-t}{2}\right)^{\frac{s}{2}} \cdot \sum_{m=0}^{N_t} \widehat{C}_{m,k} P_m^{\frac{s}{2},\frac{n}{2}-1}(t) \\
  &=\sum_{k=0}^{N_t}v(t_k)  \cdot \sum_{m=0}^{N_t}  \widetilde{C}_{m,k} \left(\frac{1-t}{2}\right)^{\frac{s}{2}} P_m^{\frac{s}{2},\frac{n}{2}-1}(t),  
\end{aligned}
\end{equation}
where the modified coefficients
$
  \widetilde{C}_{m,k} =\left(\frac{2}{1-t_k}\right)^{\frac{s}{2}}\widehat{C}_{m,k}= \frac{2^{\frac{s}{2}}}{\gamma_m^{\frac{s}{2},\frac{n}{2}-1}}(1-t_k)^{-\frac{s}{2}} \omega_k^{\frac{s}{2},\frac{n}{2}-1} P_m^{\frac{s}{2},\frac{n}{2}-1}(t_k).
$

Through the substitution $t= 2|x|^2-1$, and letting $v(t) = v(2|x|^2-1)$ in \eqref{eq:ninterp}. For convenience, we denote $u(x) = v(2|x|^2-1)$ and $u(|x_k|) = v(2|x_k|^2-1)$ where $t_k = 2|x_k|^2-1$. 
Then we have the radial spectral approximation:
\begin{equation}\label{eq:26}
  \mathcal{I}_{N_t}^{\frac{s}{2}}u(x) = \sum_{k=0}^{N_t}u(|x_k|)\sum_{m=0}^{N_t} \widetilde{C}_{m,k} (1-|x|^2)^{\frac{s}{2}}P_m^{\frac{s}{2},\frac{n}{2}-1}(2|x|^2-1), \quad x \in B_1^0.
\end{equation}
Here, $u(|x_k|)$ implies that we only need to select points $\left\{ x_k \right\}_{k=0}^{N_t}$ at a random angular direction.
Combining with \eqref{eq:eigfun}, we get the fractional Laplacian action as
  \begin{align}
    (-\Delta)^{\frac{s}{2}} \left[ \mathcal{I}_{N_t}^{\frac{s}{2}}u(x) \right] 
    &= \sum_{k=0}^{N_t}u(|x_k|) \sum_{m=0}^{N_t}\widetilde{C}_{m,k} \, \mu_m P_m^{\frac{s}{2},\frac{n}{2}-1}(2|x|^2-1) \notag\\
    &= \sum_{m=0}^{N_t} \hat{u}_m P_m^{\frac{s}{2},\frac{n}{2}-1}(2|x|^2-1),  \text{ where } \hat{u}_m = \mu_m \sum_{k=0}^{N_t} \widetilde{C}_{m,k} u(x_k).\label{eq:fracInt}
  \end{align}

According to \eqref{eq:26} and \eqref{eq:fracInt}, we can obtain the spectral Monte Carlo method in high-dimensional ball $B_1^0$. The implementation process is similar to Algorithm \ref{alg:SMC}, where the main difference is that the points to be calculated and the key formulas \eqref{eq:26} and \eqref{eq:fracInt}. For the Jacobi-Gauss points $\{ t_k \}_{k=0}^{N_t}$, then we get $|x_k| = \sqrt{(t_k+1)/2}$ and plus a random angular direction to generate $\{x_k\}_{k=0}^{N_t}$. Taking $n=2$ as an example, we employ polar coordinates and choose a random variable $\theta \in (0,2\pi)$, then the point $x_k=(r_k\cos\theta, r_k\sin\theta)$

\subsection{Comparisons, Explanations, and Remarks}\label{sec:Compa}

\begin{enumerate}
\item
  As established in \cite{GRUBB2015Fractional}, the solution $u$ of \eqref{eq:fracPoiss} admits a decomposition into a boundary singular term and a smooth component. Specifically, for a smooth source term,
\begin{equation}\label{eq:reg}
f(x) \in C^{\infty}(\bar{\Omega}) \Longleftrightarrow u \in \mathrm{dist}(x,\partial\Omega)^{s/2} C^{\infty}(\bar{\Omega}), \quad \Omega \subset \mathbb{R}^{n},
\end{equation}
where $\mathrm{dist}(x, \partial\Omega)$ denotes the Euclidean distance from $x \in \Omega$ to the boundary. This characterization precisely captures the limited regularity of the solution near $\partial\Omega$.

In the one-dimensional case with $\Omega=(-1,1)$ discussed in Section \ref{sec:1eig}, the weight function $\rho(t) = (1-t^2)^{s/2}$ in \eqref{eq:eigfun1} exactly matches the boundary singularity described in \eqref{eq:reg}. Furthermore, since the Jacobi polynomials $P_m^{s/2,s/2}(t)$ are smooth, the generalized Jacobi polynomials $J^{s/2}_m(t)$ inherently satisfy the decomposition structure:
\begin{equation}\label{eq:bound}
J^{s/2}_m(t) \in \mathrm{dist}(t,\partial\Omega)^{s/2} C^{\infty}(\bar{\Omega}) \Longleftrightarrow (-\Delta)^{s/2}J^{s/2}_m(t) = \lambda_m P_m^{s/2,s/2}(t) \in C^{\infty}(\bar{\Omega}).
\end{equation}

A key distinction is that \eqref{eq:reg} holds in arbitrary dimensions, while \eqref{eq:bound} and the associated eigenfunction relation \eqref{eq:eigfun1} are intrinsically one-dimensional. Consequently, to implement the spectral Monte Carlo method in higher dimensions ($n \geq 2$), a new eigenfunction relation analogous to \eqref{eq:eigfun1} must be developed for radial solutions.

\item
In Section \ref{sec:neig}, we derived a new eigenfunction relation \eqref{eq:eigfun}. Its weight function $\rho(x) = (1-|x|^2)^{s/2}$ for $x \in \Omega = B_1^0$ correctly captures the boundary singularity $\mathrm{dist}(x, \partial\Omega)^{s/2}$. It remains to verify the smoothness of the Jacobi polynomials $P_m^{s/2,  n/2-1}(2|x|^2-1)$.

Let $h(x) = |x|$. Its partial derivatives are $\partial h(x)/\partial x_i = x_i / |x|$ for $x = (x_1, x_2, \cdots, x_n)$. Since $1/|x|$ is not differentiable at the origin, any function involving odd powers of $|x|$ may lack smoothness there. However, the polynomial $P_m^{s/2, n/2-1}(2|x|^2-1)$ is an even function of $|x|$, as its argument depends on $|x|^2$. This even symmetry ensures that the singularities in the derivatives cancel, guaranteeing that the function is smooth on $\bar{\Omega}$.
Therefore, both the weight $\rho(x)$ and the polynomial $P_m^{s/2, n/2-1}(2|x|^2-1)$ conform to the regularity structure in \eqref{eq:reg}, satisfying:
\begin{align}
\rho(x)P_m^{s/2, n/2-1}(2|x|^2-1) & \in \mathrm{dist}(x, \partial\Omega)^{s/2} C^{\infty}(\bar{\Omega}) \notag \\
& \Longleftrightarrow (-\Delta)^{s/2} \left[\rho(x)P_m^{s/2, n/2-1}(2|x|^2-1)\right] \in C^{\infty}(\bar{\Omega}).\label{eq:reg2}
\end{align}

\item
We also note a connection established in \cite{Olver2010NIST}:
\begin{equation}
\frac{P_{2m}^{\frac{s}{2},\frac{s}{2}}(t)}{P_{2m}^{\frac{s}{2},\frac{s}{2}}(1)} = \frac{P_m^{\frac{s}{2},\frac{1}{2}-1}(2t^2-1)}{P_m^{\frac{s}{2},\frac{1}{2}-1}(1)},
\end{equation}
with the corresponding eigenvalues satisfying $\lambda_{2m} = \mu_m$. This identity demonstrates that the $n$-dimensional eigenfunction relation \eqref{eq:eigfun} for the special case $n=1$ reduces to the one-dimensional result \eqref{eq:eigfun1}, but specifically for the even-indexed eigenfunctions.

Furthermore, an equivalent formulation to \eqref{eq:eigfun1} is presented in Corollary 3.15 of \cite{Acosta2018Regular}, expressed in terms of the Gegenbauer polynomials $C_m^{\frac{s}{2}+\frac{1}{2}}(t)$. The two polynomial bases are related by the identity
\begin{equation}
P_m^{\frac{s}{2},\frac{s}{2}}(t) = \frac{(\frac{s}{2}+1)_m}{(s+1)_m}C_m^{\frac{s}{2}+\frac{1}{2}}(t),
\end{equation}
where $(a)_m = a(a+1) \cdots (a+m-1)$ denotes the Pochhammer symbol.

\end{enumerate}

\section{Numerical experiments}\label{sec:numexp}
In this section, we provide some numerical examples to illustrate the effectiveness of the spectral Monte Carlo method for solving fractional Poisson problem  in $n$-dimensional space. Let $\Omega = B^0_1$ and $g(x)=0$ on $\Omega^c$. Define the error   
%\begin{equation}
$ E_{N_t}^{\infty} = \max_{0\leq k\leq N_t} \left| \mathbb{E}\left(u(x_k) - u_{\star}(x_k) \right)\right|$,
%\end{equation}
where the points $\{ x_k\}_{k=0}^{N_t}$ matches the Jacobi-Gauss quadrature nodes $\{ t_k \}_{k=0}^{N_t}$ with index $\left(\frac{s}{2},\frac{n}{2}-1\right)$ and $ t_k = 2|x_k|^2 -1$.  For convenience, we explain the meaning of some notations as follows:
\begin{itemize}
  \item $N_t$: the degree of Jacobi polynomials.
  \item $M$: the number of Monte Carlo trials.
  \item $K$: the iteration number.
\end{itemize}

\begin{exm}
Let $n=1$ so $\Omega = (-1,1)$. 
According to characteristic relation \eqref{eq:eigfun}, if we take the source function 
\begin{equation*}
f(x)= (-\Delta)^{\frac{s}{2}}(1-x^2)^{\frac{s}{2}} (x^2+1) = \left(\mu_1 P_{1}^{\frac{s}{2},-\frac{1}{2}}(2x^{2}-1) + \frac{1}{2}\mu_0 P_{0}^{\frac{s}{2},-\frac{1}{2}}(2x^{2}-1)\right) \Big/ \left( \frac{s+3}{2} \right).
\end{equation*}  
the true solution of the equation is $u(x) = (1-x^2)^{\frac{s}{2}} (x^2+1)$.
\end{exm}

\begin{figure}[htpb]
  \centering
  \subfloat[$N_t=2$, $M = 50$.]{
      \label{fig:2a}
      \includegraphics[width=0.35\textwidth]{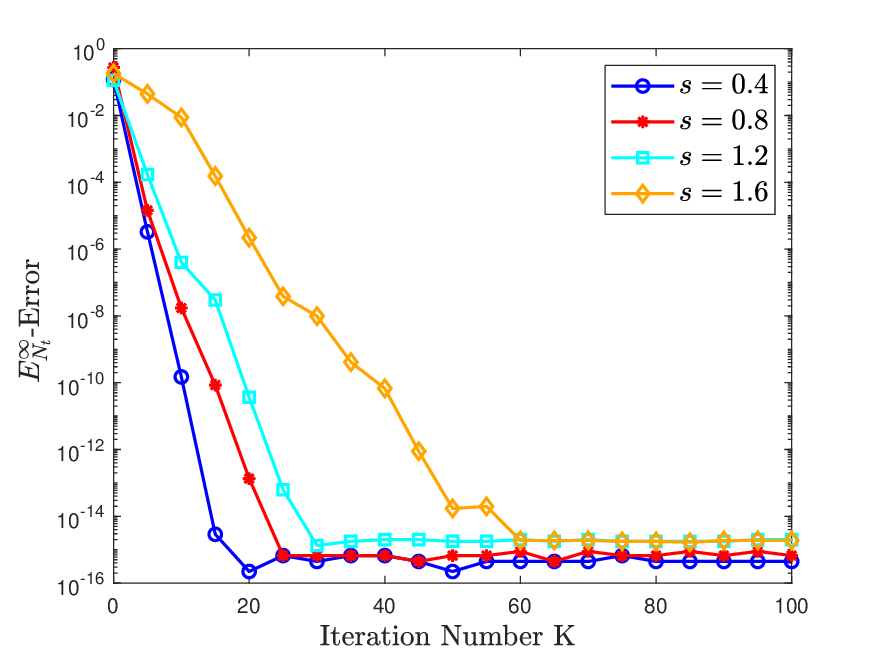}}
  \subfloat[$N_t=2$, $s = 0.4$.]{
      \label{fig:2b}
      \includegraphics[width=0.35\textwidth]{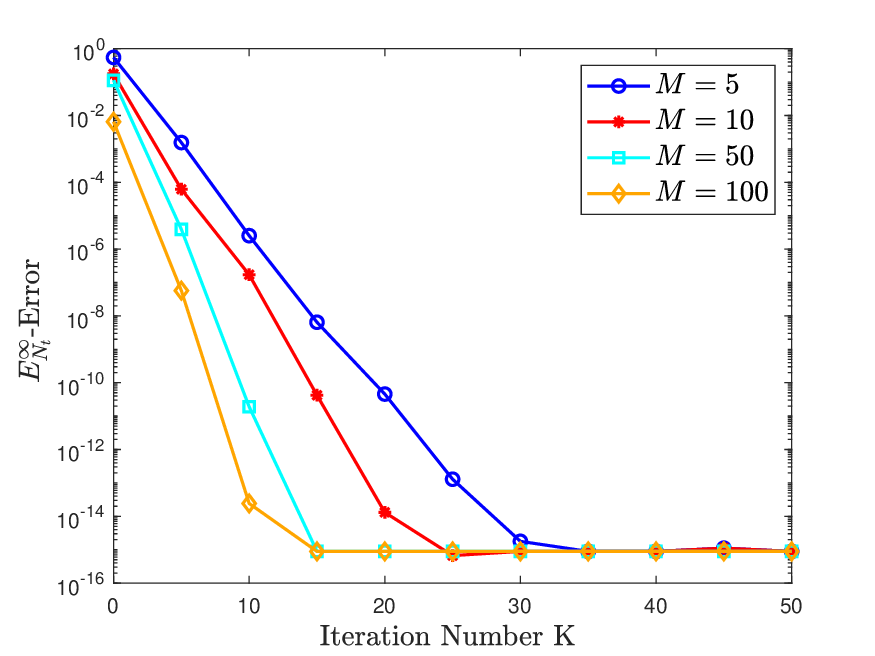}}
  \caption{$n=1$, the discrete maximum error $E_{N_t}^{\infty}$ about iterative step number $K$.}
\label{fig:2}
\end{figure}

The purpose of this example is to demonstrate the effectiveness of our method for $1$ dimension equation. In Figure \ref{fig:2a}, we fix $M = 50$ and $N_x=2$. For the various fractional parameters $s=0.4,0.8,1.2,1.6$, we observe that the error $E_{N_t}^{\infty}$ for a given $s$ decays exponentially as iteration number $K$ increases, and reaches the spectral accuracy more quickly when $s$ is smaller. 
In Figure \ref{fig:2b}, we choose to fix $N_x$ and $s=0.4$, and the exponential decay is also found. As the number of Monte Carlo experiments $M$ increases, the error $E_{N_t}^{\infty}$ achieves spectral accuracy faster, which accord with the feature of Monte Carlo method.

\begin{exm}\label{exm:2}
Let $n=2$ and $f(x) = 2^{s}\Gamma(1+\frac{s}{2})\Gamma(\frac{n+s}{2})/\Gamma(\frac{n}{2})$. Then we have $u(x) = (1-|x|^2)^{\frac{s}{2}}$; See for example \cite{Dyda2012power}.
\end{exm}

\begin{figure}[htpb]
  \centering
  \subfloat[$N_t=2$, $M = 50$.]{
      \label{fig:4a}
      \includegraphics[width=0.35\textwidth]{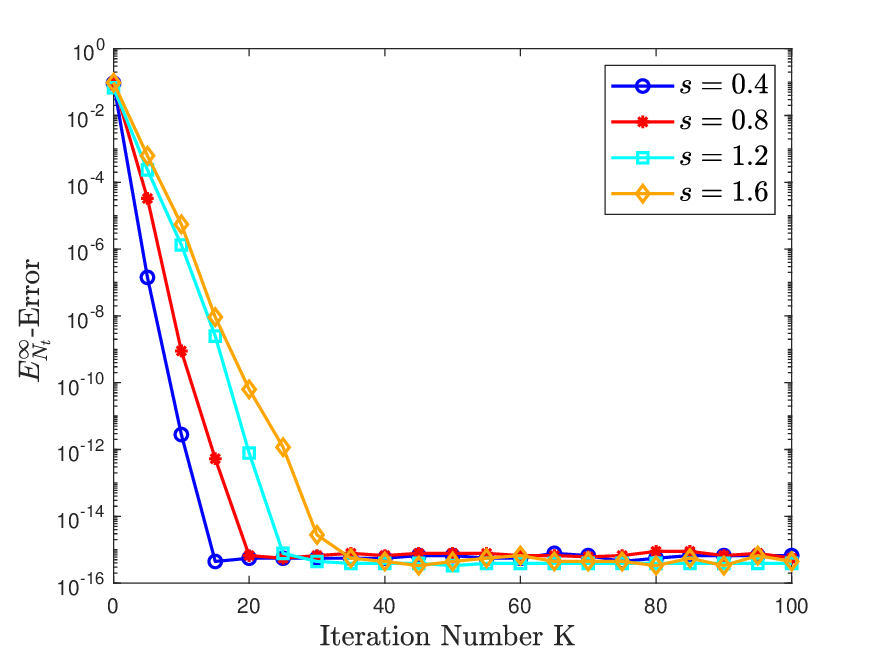}}
  \subfloat[$N_t=2$, $s = 0.4$.]{
      \label{fig:4b}
      \includegraphics[width=0.35\textwidth]{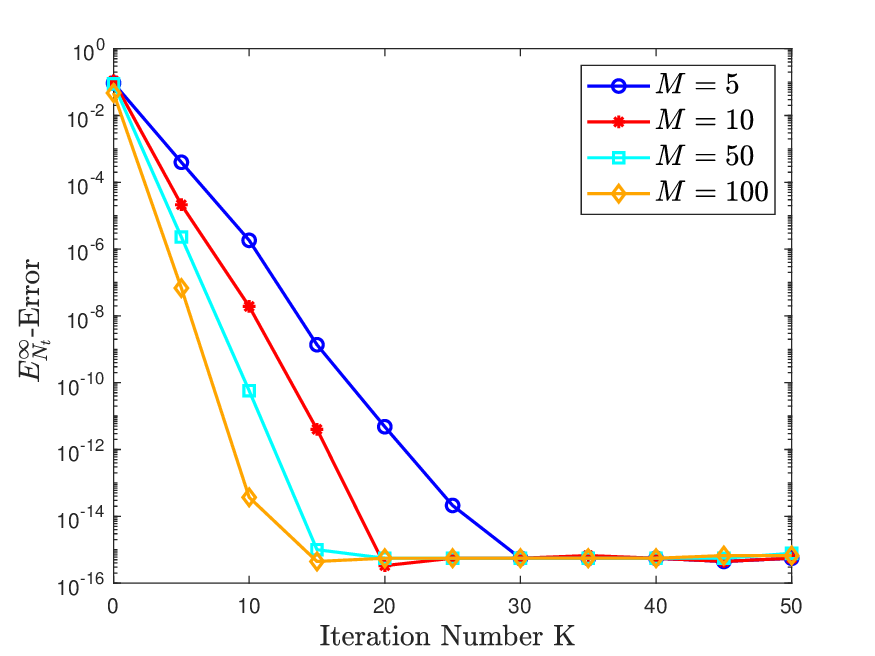}}
  \caption{$n=2$, the discrete maximum error $E_{N_t}^{\infty}$ for Example \ref{exm:2} about iterative steps $K$.}
  \label{fig:4}
\end{figure}

\begin{figure}[H]
  \centering
  \subfloat[$s = 0.4$]{
      \includegraphics[width=0.24\textwidth]{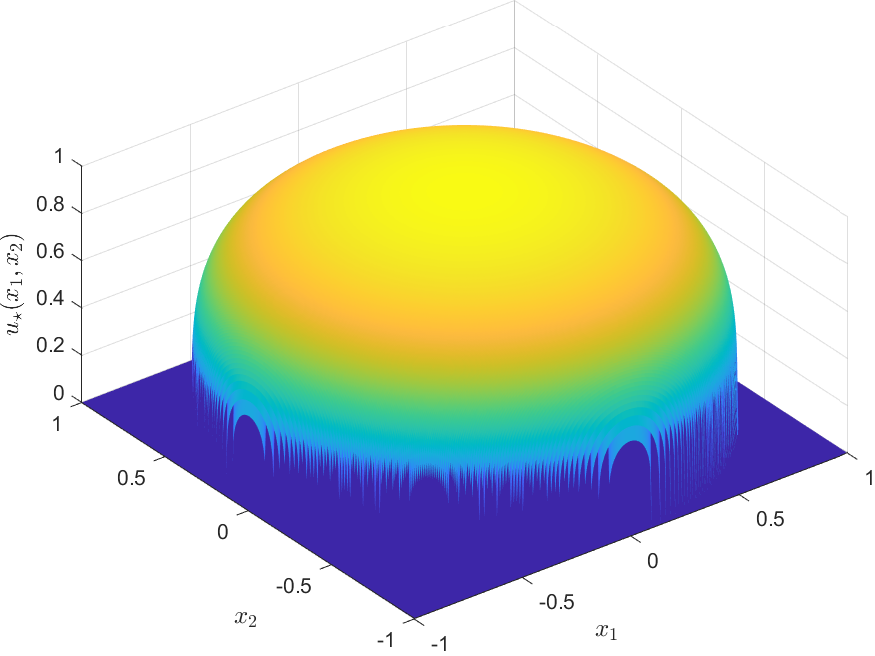}}
  \subfloat[$s = 0.8$]{
      \includegraphics[width=0.24\textwidth]{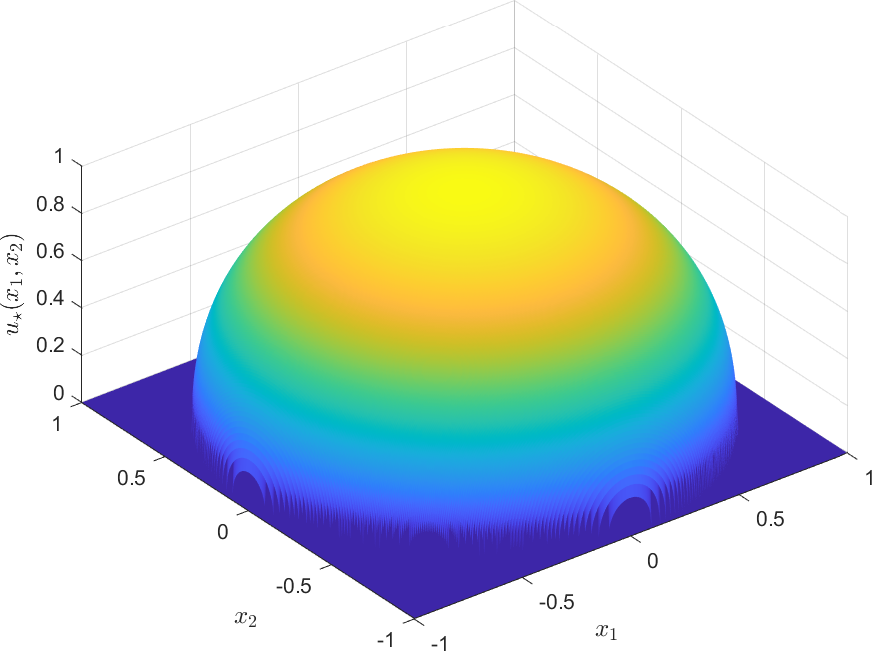}} 
  \subfloat[$s = 1.2$]{
      \includegraphics[width=0.24\textwidth]{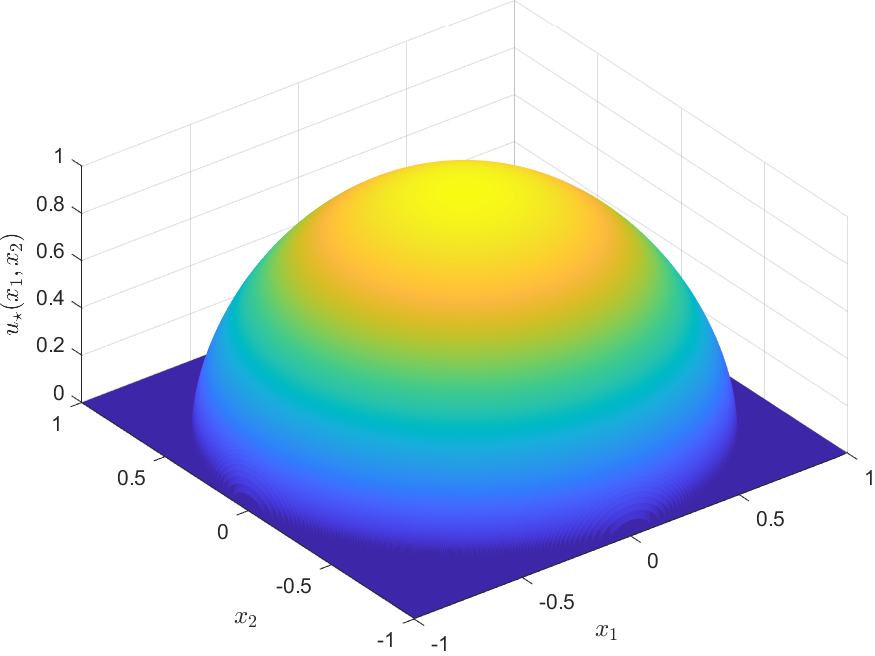}}
  \subfloat[$s = 1.6$]{
      \includegraphics[width=0.24\textwidth]{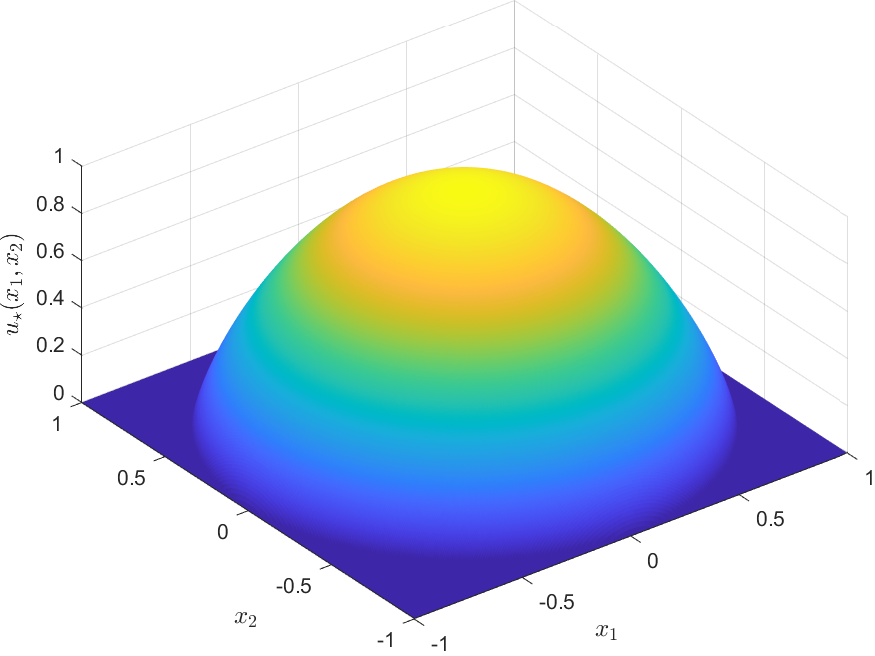}}
  \caption{$n=2$, the numerical solution $u_{\star}(x)$ for Example \ref{exm:2} with various fractional power $s$.}
  \label{fig:5}
\end{figure}

For this example, we choose a random angular direction to fix $\left\{ x_k \right\}_{k=0}^{N_t}$ on the radial direction. For the Figure \ref{fig:4a} and \ref{fig:4b}, we still observe the phenomenon of exponential convergence. Moreover, for the given $N_t$ and $M$ in Figure \ref{fig:4a}, we can use fewer iteration number $K$ to achieve spectral accuracy when $s$ is smaller. For the fixed $N_t$ and $s$, Figure \ref{fig:4b} implies that we can enhance the speed of reaching spectral accuracy by increasing the number of Monte Carlo trials $M$.
In order to verify whether the interpolation formula \eqref{eq:26} can portray the exact solutions, we exhibit the image of numerical solutions $u_{\star}(x)$ for $s = 0.4,0.8,1.2,1.6$ in the Figure \ref{fig:5}with the parameters $M=50$, $K=100$ and $N_x=2$. The singularity at boundary $\partial\Omega$ for the exact solution behaves $\mathrm{dist}(x,\partial\Omega)^{\frac{s}{2}}$, it implies that the boundary singular layer will be thinner as the fractional power $s$ becomes smaller. Our numerical results reflects this feature, i.e., we can see the sharp shape when the $s$ is smaller. This is similar to the numerical results observed in \cite{Sheng2023Efficient}.

\begin{exm}\label{exm:3}
Let $n=2$, and we choose $f(x) = \sin(|x|^2)$. The exact solution is unknown in this example.
\end{exm}

\begin{figure}[htpb]
  \centering
  \subfloat[$N_t=12$, $M = 2000$.]{
      \label{fig:8a}
      \includegraphics[width=0.35\textwidth]{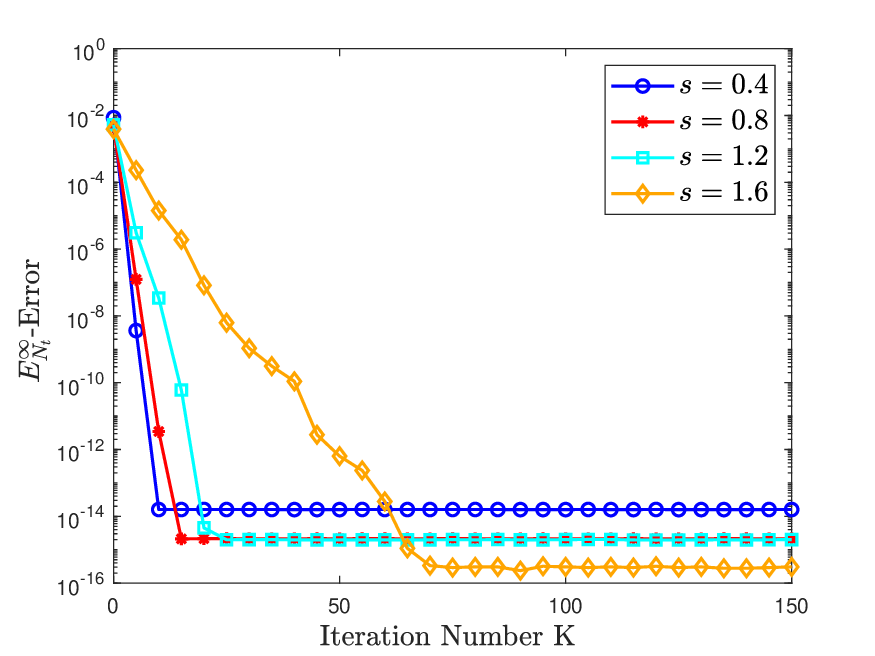}}
  \subfloat[$N_t=10$, $s = 0.4$.]{
      \label{fig:8b}
      \includegraphics[width=0.35\textwidth]{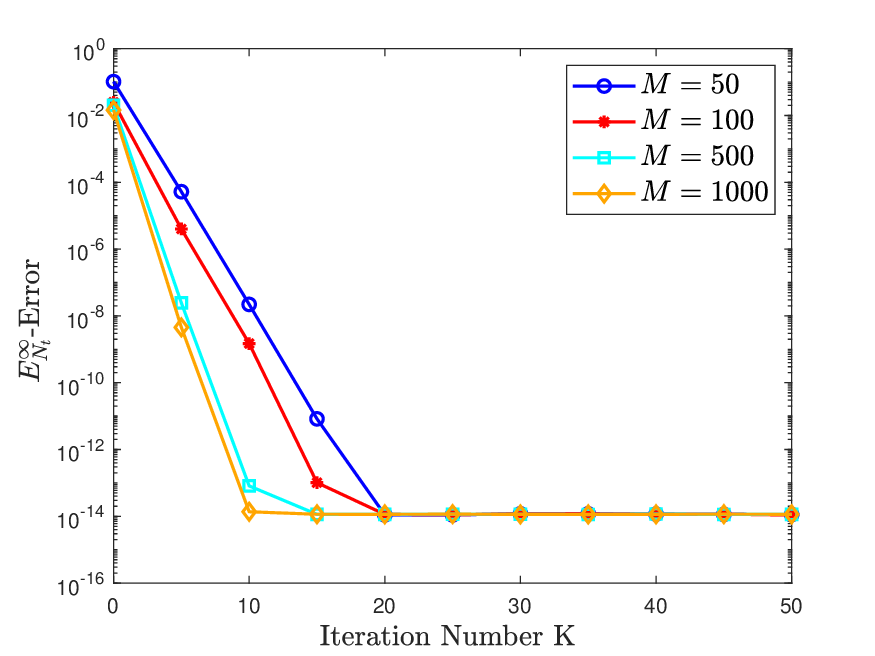}}
  \caption{$n=2$, the discrete maximum error $E_{N_t}^{\infty}$ for Example \ref{exm:3} about iterative steps $K$.}
  \label{fig:8}
\end{figure}

\begin{figure}[H]
  \centering
  \subfloat[$s = 0.4$, $N_t=12$, $M=50$]{
      \includegraphics[width=0.24\textwidth]{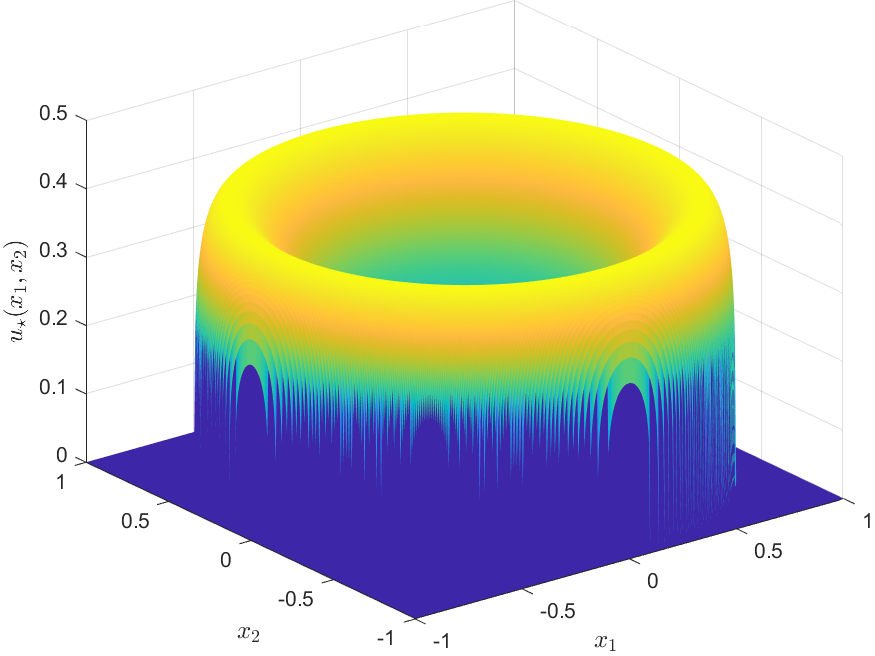}}
  \subfloat[$s = 0.8$, $N_t=12$, $M=50$]{
      \includegraphics[width=0.24\textwidth]{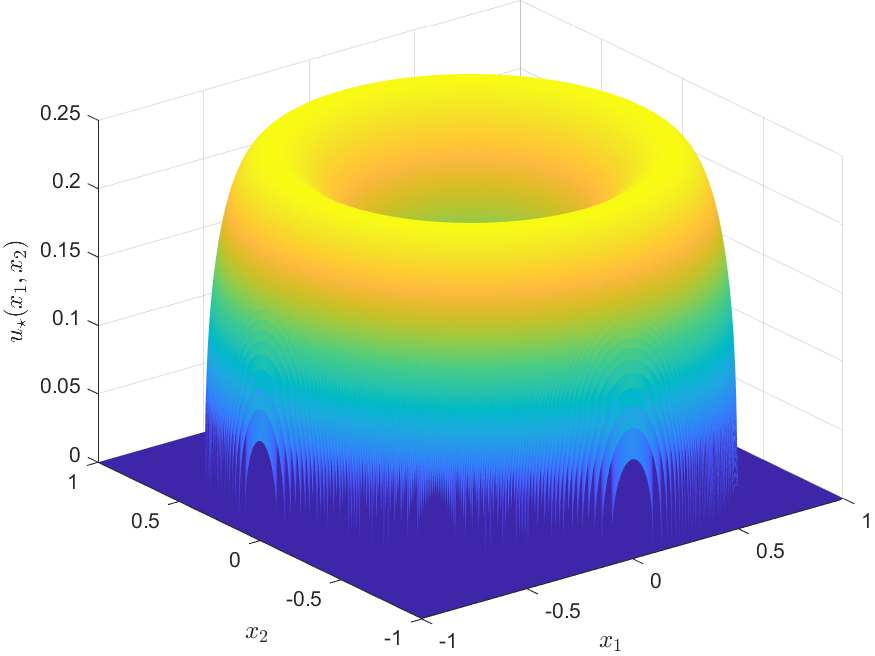}}
  \subfloat[$s = 1.2$, $N_t=12$, $M=200$]{
      \includegraphics[width=0.24\textwidth]{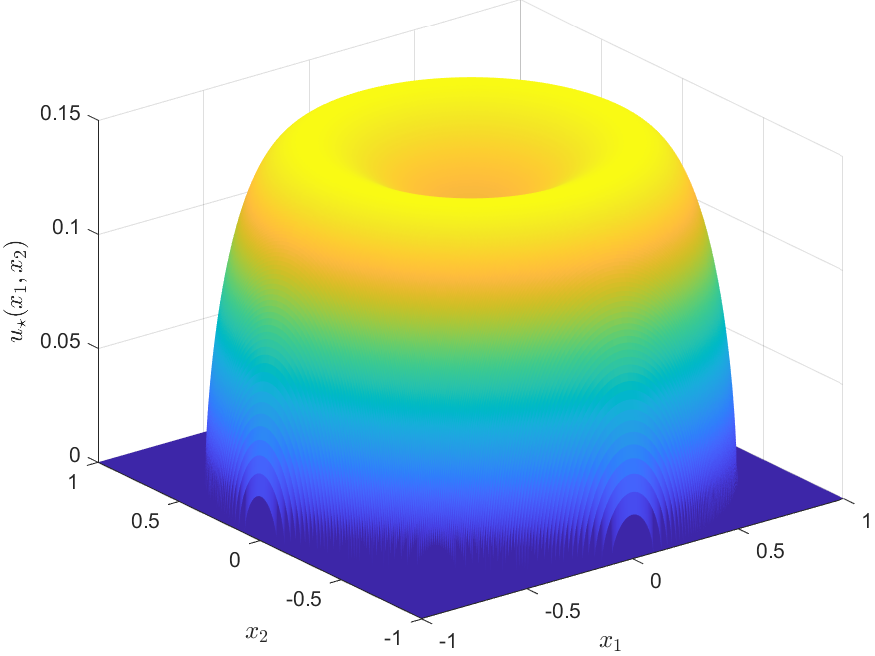}}
  \subfloat[$s = 1.6$, $N_t=12$, $M=2000$]{
      \includegraphics[width=0.24\textwidth]{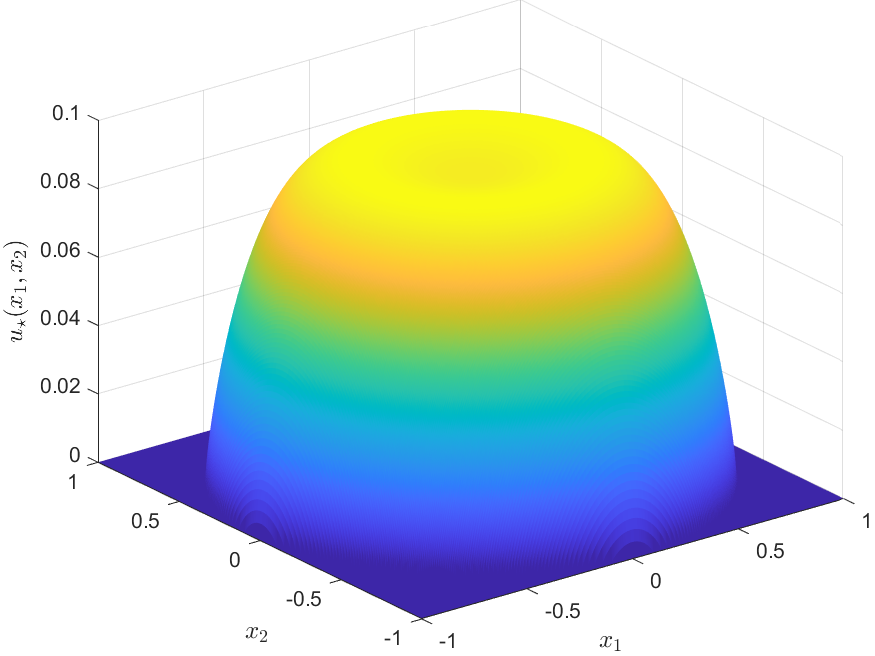}}
  \caption{$n=2$, the numerical solution $u_{\star}(x)$ for Example \ref{exm:3} with various fractional power $s$.}
  \label{fig:6}
\end{figure}

Since the exact solution for example \ref{exm:3} is unknown, we employ the spectral Galerkin method stated in \cite{hao2021sharp} to obtain the reference solution. When the solution is restricted as radial function, this spectral Galerkin method can be simplified as the method provided in \cite{xu2018spectral}. 

In Figure \ref{fig:8a}, the exponential convergence is also observed, and we can obtain high accuracy numerical solution. However, it still exists some difference with the examples containing exact solution. The first one is that the value of steady state error $E_{N_t}^{\infty}$ for various $s$ is inconsistent. But we need to point that this difference is not our main concern, it is important that the error $E_{N_t}^{\infty}$ can reach stablity quickly as $s$ decreases, this feature is same as other examples. 
The second one is that in the specific experiments, we only take the number of Monte Carlo paths $M=50$ or $100$ to make the numerical solution converge when $s$ is small $(s<1)$. For the larger $s$ $(s>1)$, we usually need a much bigger $M$ to obtain the convergence as $s$ increases. For example, $M=200$ when $s=1.2$ and $M=1000$ when $s=1.6$.  
The reason of this phenomenon is uncertain, but we guess that it may be related to hypersingularity of fractional Laplacian operator when $s>1$.
In Figure \ref{fig:8b}, we observe the same behavior as other examples that we can use fewer iteration number $K$ to get the spectral accuracy by adding $M$.

We also plot the numerical solutions $u_{\star}(x)$ for $s = 0.4,0.8,1.2,1.6$ in the Figure \ref{fig:6}. The other parameters are fixed as $M=2000$, $K=100$ and $N_x=12$. Regarding the boundary singular layer determined by $\mathrm{dist}(x,\partial\Omega)^{\frac{s}{2}}$, we also see the sharp shape when the $s$ is smaller, i.e., the boundary singular layer will become thinner as the fractional power $s$ decreases.

\begin{exm}\label{exm:4}
Let $n=10$ and  $f(x) = 2^{s}\Gamma\left(2+\frac{s}{2}\right)\Gamma\left(\frac{n+s}{2}\right)\Gamma\left(\frac{n}{2}\right)^{-1}\left(1-\left(1+\frac{s}{n}\right)|x|^2\right)$. Then we have  $u(x) = (1-|x|^2)^{1+\frac{s}{2}}$ from \cite{Dyda2012power}.
\end{exm}

\begin{figure}[htpb]
  \centering
  \subfloat[$N_t=2$, $M = 50$.]{
      \label{fig:7a}
      \includegraphics[width=0.35\textwidth]{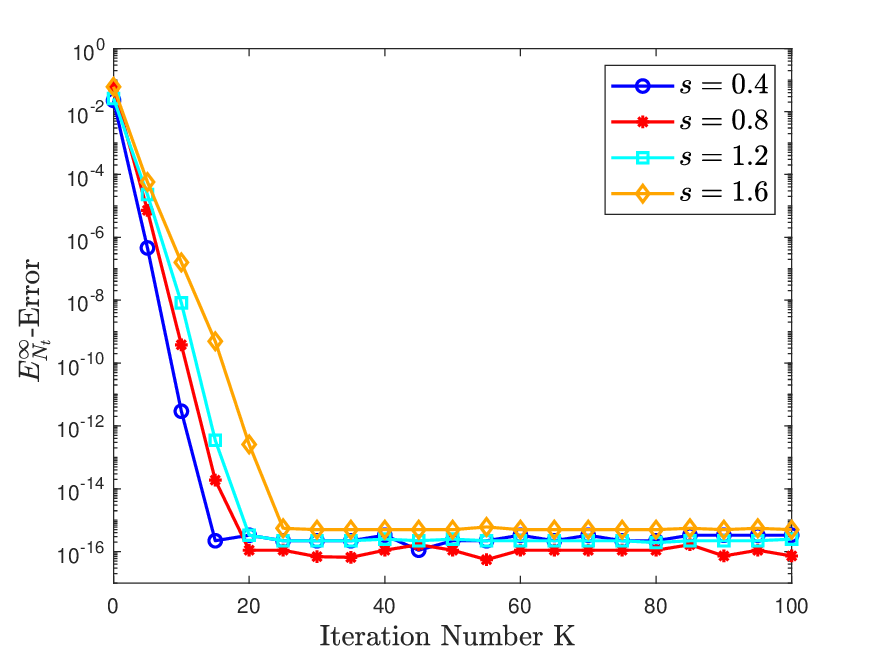}}
  \subfloat[$N_t=2$, $s = 0.4$.]{
      \label{fig:7b}
      \includegraphics[width=0.35\textwidth]{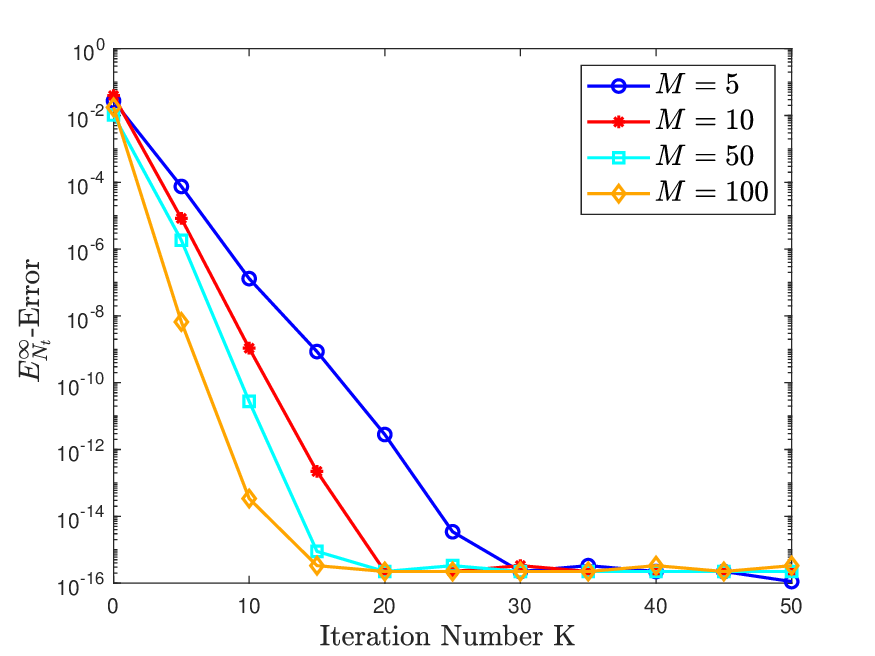}}
  \caption{$n=10$,  the discrete maximum error $E_{N_t}^{\infty}$ for Example \ref{exm:4} about iterative steps $K$.}
  \label{fig:7}
\end{figure}

In order to test the effectiveness of our method in high-dimensional space, we select a example with $n=10$.
Similarly, we also observe that the errors $E_{N_t}^{\infty}$ can maintain exponential decay from Figure \ref{fig:7a} and \ref{fig:7b}. This two picture also shows that we can reach the spectral accuracy quickly with small $s$ or large Monte Carlo trial number $M$. Generally speaking, the numerical results can well illustrate that our strategy possess the capability to carry out in high-dimensional situations. 

\section{Concluding remarks}\label{sec:conclude}
The spectral Monte Carlo method for fractional Possion problem in one dimensional spaces was first developed in \cite{Feng2025ExponentiallyAS}. Inspired by this work and another radial eigenfunction relation, we provide a new interpolation formula to make spectral Monte Carlo method suitable for high-dimensional ball region. Numerical results demonstrate our method is effective.

\bibliographystyle{alpha}
\bibliography{reference}

\end{document}